\documentclass{amsart}
\usepackage{graphicx,xcolor,textpos}
\usepackage{helvet}
\usepackage{graphicx}
\usepackage{epstopdf}
\usepackage{tikz-cd}
\usepackage{amsmath}
\usepackage{amssymb}
\usepackage{amsfonts}
\usepackage{amsthm}
\usepackage{hyperref}
\usepackage{dsfont}
\usepackage{bbm}
\usepackage{float}
\usepackage{caption}
\usepackage{subcaption}
\usepackage{enumitem}	
\usepackage{verbatim}
\usepackage[capitalise]{cleveref}

\usepackage[ruled,vlined]{algorithm2e}
\usepackage[utf8]{inputenc}
\usepackage[intoc]{nomencl}

\newtheorem{thm}{Theorem}[section]

\newtheorem{lem}{Lemma}[section]

\theoremstyle{definition}

\title[Tensor completion using geodesics on Segre manifolds]{Tensor completion using geodesics on Segre manifolds}

\author[Lars Swijsen]{Lars Swijsen}
\address[Lars Swijsen]{KU Leuven, Department of Mathematics, Celestijnenlaan 200B -- Box 2400, B-3001 Leuven, Belgium}
\email[Lars Swijsen]{lars.swijsen@kuleuven.be}

\author[Joeri Van der Veken]{Joeri Van der Veken}
\address[Joeri Van der Veken]{KU Leuven, Department of Mathematics, Celestijnenlaan 200B -- Box 2400, B-3001 Leuven, Belgium}
\email[Joeri Van der Veken]{joeri.vanderveken@kuleuven.be}
\thanks{J. V. d. V. is supported by the Research Foundation -- Flanders (FWO) and the Fonds de la Recherche Scientifique (FNRS) under EOS Project G0H4518N}

\author[Nick Vannieuwenhoven]{Nick Vannieuwenhoven}
\address[Nick Vannieuwenhoven]{KU Leuven, Department of Computer Science, Celestijnenlaan 200A -- Box 2402, B-3001 Leuven, Belgium}
\email[Nick Vannieuwenhoven]{nick.vannieuwenhoven@kuleuven.be}
\thanks{N. V.~was supported by a Postdoctoral Fellowship of the Research Foundation---Flanders (FWO) with project 12E8119N. This work was initiated while N. V. visited Carlos Beltr\'an (Universidad de Cantabria), which was supported by the FWO Grant for a long stay abroad V401518N}


\begin{document}
\begin{abstract}
We propose a Riemannian conjugate gradient (CG) optimization method for finding low rank approximations of incomplete tensors. Our main contribution consists of an explicit expression of the geodesics on the Segre manifold. These are exploited in our algorithm to perform the retractions. We apply our method to movie rating predictions in a recommender system for the MovieLens dataset, and identification of pure fluorophores via fluorescent spectroscopy with missing data. In this last application, we recover the tensor decomposition from less than $10\%$ of the data.
\end{abstract}
\maketitle

\section{Introduction}

Data is often stored in a $d$-dimensional array $T = [T_{i_1,\dots,i_d}]_{i_1=1,\dots,i_d=1}^{n_1,\dots,n_d}$, called a \textit{tensor}. Let $a_1,\dots,a_d$ be vectors with $a_i\in \mathbb{R}^{n_i}$. We can turn these vectors into a tensor by applying the outer product and defining \[T_{i_1,\dots,i_d} = a_1(i_1)\cdot\dots \cdot a_d(i_d).\]We call this tensor a \emph{rank-$1$ tensor} and denote it by $a_1\otimes \dots\otimes a_d$. 

Because of the high dimensionality of a tensor, one is often interested in storing the data in a concise way. A classic way to achieve this is using a \textit{canonical polyadic decomposition} (CPD). Given a tensor $T$ it is possible to write it as a finite sum of rank-$1$ tensors \[ T=\sum_{j=1}^r a_1^j\otimes \dots\otimes a_d^j \] with $a_i^j\in \mathbb{R}^{n_i}$. The minimal number $r$ of rank-$1$ tensors needed to build the tensor $T$ is called the \textit{rank} of $T$. 

The CPD has applications in different fields, such as signal processing, chemometrics and linguistics \cite{kolda_tensor_2009,bro_parafac_1997}. It turns out that because of extra underlying structure in these tensors, the rank-$1$ tensors that are obtained contain valuable information. Using the fact that for low ranks the rank-$1$ terms are essentially unique \cite{kruskal_three-way_1977,sidiropoulos_uniqueness_2000,chiantini_algorithm_2014}, a property that is called \textit{r-identifiability}, one can recover this information from the CPD.

In this paper we present a Riemannian optimization algorithm that seeks a best rank-$r$ approximation of a given tensor $T$. Concretely, we aim to solve the following optimization problem:  \[ \text{Given an arbitrary tensor }T, \]
\[ \text{find rank-$1$ tensors } T_1,\dots, T_r \text{ such that} \]\[
\left\|T-\sum_{j=1}^rT_j\right\|_\text{F} \text{ is minimal}, \]
where $\left\|T\right\|_\text{F} = \sqrt{\sum_{i_1=1,\dots,i_d=1}^{n_1,\dots,n_d} T_{i_1,\dots,i_d}^2}$  is called the \textit{Frobenius norm} of $T$ \footnote{We assume a solution for this optimization problem exists. However, this is not always the case, as is explained in \cite{de_silva_tensor_2008}}. Therefore, we optimize over $r$ copies of the space of all rank-$1$ tensors. It turns out that the space of all tensors of rank $1$ has an interesting structure. It is a space with curvature, and more specifically it is a \textit{Riemannian manifold} called the \textit{Segre manifold}, which we denote by $Seg$. In this paper we compute the geodesics of this Riemannian manifold and obtain the following result:
\begin{thm}
\label{geo}
The (unit-speed) geodesic of $Seg$ through the point \[p=(\lambda,x_1,\dots,x_d) \in \mathbb{R}^+_0\times \mathbb{S}^{n_1-1}\times \dots \times \mathbb{S}^{n_d-1}\] in the direction \[(\dot \lambda,\dot x_1,\dots,\dot x_d) \in T_p\left(\mathbb{R}^+_0\times \mathbb{S}^{n_1-1}\times \dots \times \mathbb{S}^{n_d-1}\right)\] is given by
\begin{align*}
\lambda(t)&= \sqrt{t^2+\frac{2\lambda P}{\sqrt{P^2+1}}t+\frac{\lambda^2}{P^2+1}}, \\
x_i(t) &= x_i\cos\left(\frac{||\dot{x_i}||}{M}f(t)\right)
+\frac{\dot x_i}{||\dot x_i||}\sin\left(\frac{||\dot{x_i}||}{M}f(t)\right), \qquad i=1,\dots,d,
\end{align*}
where $f(t)=\tan^{-1}\left(\frac{\sqrt{P
^2+1}}{\lambda}t+P\right)-\tan^{-1}(P)$, $P=\frac{\dot\lambda}{\lambda M}$ and $M=\sqrt{\sum_{k=1}^d||\dot x_k||^2}$.
\end{thm}
These geodesics, while interesting in their own right, are then used to implement a Riemannian conjugate gradient method that moves over the manifold in a geometrically optimal way. These geodesics are a computationally viable alternative to HOSVD-based retractions that can be found in the literature, see for example \cite{breiding_riemannian_2018,kressner_low-rank_2014}.

As an application of these optimal curves, we investigate two applications involving missing data, i.e.,  \textit{incomplete} tensors. First, we will look at the 1M MovieLens database. This is a collection of approximately 1 million ratings given by 6000 users on 4000 movies. We aim to use the algorithm to create a prototype recommender system that predicts the match between an arbitrary user-movie pair, as in \cite{frolov_tensor_2017}. Next, we apply the algorithm in the setting of fluorescence spectroscopy. To measure fluorescent compounds in a solution, the latter is excited with light of different wavelengths. This excites the electrons and will result in the emission of light of different wavelengths. The resulting emission patterns can be used to determine the fluorescent compounds in the mixtures. We experimentally study how incomplete the data can become, such that it remains possible to recover the entire dataset. These results might be applied in developing new compressive sensing hardware that reduces the experimental duration.

The paper is structured as follows. In the next section we give an short introduction to Riemannian geometry. In particular, the relevant definitions and notation will be introduced. \Cref{sectiongeo} contains the computation of the geodesics of the Segre manifold. The geodesics are used in \cref{sectionalgo} to create a Riemannian conjugate gradient algorithm for tensor decomposition, whose experimental results are discussed in \cref{sectionresults}.

\subsection*{Acknowledgements}
We thank Carlos Beltr\'an from the Universidad de Cantabria for initiating the computation of the geodesics of the Segre manifold while the last author visited him in the first semester of $2018$. 

\section{Introduction to Riemannian Geometry}
\label{introductionRiemann}

In this section we will give a brief introduction to some of the elementary notions of Riemannian geometry that we need in order to present the main result. If the reader is interested in a more in-depth introduction to this subject, we refer to \cite{godinho_introduction_2014}.

A \textit{manifold} $M$ is a space that locally looks like a Euclidean space, together with some topological conditions. It is a second-countable and Hausdorff space that is \textit{locally  homeomorphic} \footnote{This means $M$ can be covered by open sets $U_i$ and continuous maps $\phi_i:U_i \to \mathbb{R}^n$, called \textit{charts}, such that for each $i$, $\phi_i$ is invertible and $\phi_i^{-1}$ is continuous. Additionally, if a $U_i$ and $U_j$ overlap, then the map $\phi_j\circ\phi_i^{-1}$ has to be differentiable on the overlap.} to an Euclidean space. If the manifold $M$ locally resembles $\mathbb{R}^n$, it is said to be of \textit{dimension} $n$. At each point $p \in M$ one can define the \textit{tangent space} of $M$ at $p$, denoted by $T_pM$. This linear space contains all velocity vectors of paths in $M$ at $p$. Concretely it is defined by an equivalence relationship on all paths $\gamma: \mathbb{R} \to M$  with $\gamma(0)=p$ where $[\gamma_1]=[\gamma_2]$ if and only if $\gamma_1'(0) = \gamma_2'(0)$.\footnote{Choose a chart $\phi_i:U_i \to \mathbb{R}^n$ such that $p\in U_i$. Then $\gamma'(0)$ is defined as the regular derivative $\frac{d}{dt}(\phi\circ\gamma)(0)$.} The union of all tangent spaces to $M$ is denoted by $TM$ and is called the tangent bundle. A \textit{vector field} $X$ on the manifold $M$ is a smooth collection of tangent vectors at every point of the manifold. It is a smooth map $X: M \to TM$ with $X(p) \in  T_pM$. Later on, we will need a tool to connect different tangent spaces on $M$ in order to transport tangent vectors from one point to another point. This process can be achieved by a \textit{connection} $\nabla$ on the manifold. We again refer to \cite{godinho_introduction_2014} for the technical definition. 

A manifold $M$ is called a \textit{Riemannian manifold} when it is endowed with a metric. A \textit{Riemannian metric} is a smooth positive-definite inner product $g_p$ on each tangent space $T_pM$. Once we have defined a Riemannian metric on $M$, there is a unique connection that interacts with this metric in a nice way. This canonical connection is called the \textit{Levi--Civita connection}. A Riemannian metric allows us to properly study notions such as distances and angles on a curved manifold. The length of a path $\gamma: I\subset \mathbb{R} \to M$ is defined as 
\[ l(\gamma) = \int_I||\gamma'(t)||_{\gamma(t)\ }\mathrm{d}t = \int_I \sqrt{g_{\gamma(t)}\left( \frac{d\gamma}{dt},\frac{d\gamma}{dt}\right)}\mathrm{d}t.\]
An important concept for the remaining sections will be the notion of a geodesic. A \textit{geodesic} is the generalisation of a straight line in Euclidean space to curved manifolds. Locally on the manifold they correspond to distance minimizing curves. In \cite{godinho_introduction_2014} one can find a definition which involves the Levi--Civita connection. 

\section{Geodesics on the Segre manifold}
\label{sectiongeo}

In this section we find the geodesics of the Segre manifold consisting of all tensors of rank $1$, equipped with the metric induced by the ambient Euclidean metric. Thus, we are looking for the length-minimizing curves in the usual Euclidean distance in the space of all tensors. For simplicity of notation, we restrict ourselves in the introduction of this section to the case of tensors of order $3$ with size $n_1\times n_2\times n_3$. The main theorem will be proved for tensors of arbitrary order and size. 

Let us define $\mathbb{S}^n$ as the unit sphere in $\mathbb{R}^{n+1}$, i.e., the manifold of points in $\mathbb{R}^{n+1}$ with Euclidean distance $1$ to the origin. Let $\mathbb{R}_0^+$ denote the manifold of strictly positive real numbers. Note that, using the multilinear properties of the tensor product, an arbitrary rank-$1$ tensor of size $n_1\times n_2 \times n_3$ can be parametrized by the following local diffeomorphism:\footnote{A map $f:M\to N$ is a local diffeomorphism if for each $p \in M$ there exists an open set $U$ around $p$ such that $f(U)$ is open in $N$ and the restriction of $f$ is differentiable with a differentiable inverse function.}
\begin{align*}
\phi: \mathbb{R}_0^+ \times \mathbb{S}^{n_1-1} \times \mathbb{S}^{n_2-1} \times \mathbb{S}^{n_3-1} &\to Seg \\
(\lambda,a_1,a_2,a_3) &\mapsto \lambda  a_1\otimes a_2 \otimes a_3 .
\end{align*} 
This allows us to locally identify $Seg$ with $\mathbb{R}_0^+ \times \mathbb{S}^{n_1-1} \times \mathbb{S}^{n_2-1} \times \mathbb{S}^{n_3-1}$. Additionally, deriving $\phi$ gives a description of the tangent space of the Segre manifold. Let $p=(\lambda,a_1,a_2,a_3)$, then we can identify the tangent space of $Seg$ with \[T_{p}(\mathbb{R}_0^+ \times \mathbb{S}^{n_1-1} \times \mathbb{S}^{n_2-1} \times \mathbb{S}^{n_3-1}) \cong \mathbb{R} \times T_{a_1}\mathbb{S}^{n_1-1}  \times T_{a_2}\mathbb{S}^{n_2-1}\times T_{a_3}\mathbb{S}^{n_3-1}\] as follows: 
\begin{equation*}
\begin{split}
\mathrm{d}_{p}\phi : \mathbb{R} \times T_{a_1}\mathbb{S}^{n_1-1}  \times T_{a_2}\mathbb{S}^{n_2-1}\times T_{a_3}\mathbb{S}^{n_3-1}
&\to T_{\lambda a_1\otimes a_2 \otimes a_3}Seg \\
(\dot\lambda,\dot x_1,\dot x_2, \dot x_3) &\mapsto \dot\lambda a_1\otimes a_2\otimes a_3 + \lambda (\dot x_1\otimes a_2\otimes a_3 \\
& \enskip\quad + a_1\otimes\dot x_2 \otimes a_3 + a_1\otimes a_2\otimes \dot x_3).
\end{split}
\end{equation*}
Note that the Segre manifold lives in the ambient Euclidean space $\mathbb{R}^{n_1\times n_2\times n_3}$ of all tensors of size $n_1\times n_2 \times n_3$. This allows us to equip the Segre manifold with the induced metric of the ambient Euclidean space. Fix a point $(\lambda,a_1,a_2,a_3) \in \mathbb{R} \times \mathbb{S}^{n_1-1} \times \mathbb{S}^{n_2-1} \times \mathbb{S}^{n_3-1}$ and take two tangent vectors $(\dot\lambda,\dot x_1, \dot x_2, \dot x_3)$ and $(\dot\mu,\dot y_1, \dot y_2, \dot y_3)$ at this point. Using the fact that $T_p\mathbb{S}^n$ is the orthogonal complement of the position vector for all $p \in \mathbb{S}^n \subset \mathbb{R}^{n+1}$, it follows that for the ambient metric:
\begin{align*}
&\langle (\dot\lambda,\dot x_1, \dot x_2, \dot x_3),  (\dot\mu,\dot y_1, \dot y_2, \dot y_3)\rangle_{(\lambda,a_1,a_2,a_3)}  \nonumber\\
&= \langle \dot\lambda a_1\otimes a_2\otimes a_3 + \lambda (\dot x_1\otimes a_2\otimes a_3 + a_1\otimes\dot x_2 \otimes a_3 + a_1\otimes a_2\otimes \dot x_3), \nonumber\\
&\enskip\quad\dot\mu a_1\otimes a_2\otimes a_3 + \lambda (\dot y_1\otimes a_2\otimes a_3 
+ a_1\otimes\dot y_2 \otimes a_3 + a_1\otimes a_2\otimes \dot y_3) \rangle \nonumber\\
&= \langle \dot\lambda a_1\otimes a_2\otimes a_3,\dot\mu a_1\otimes a_2\otimes a_3 \rangle + \lambda^2\langle \dot x_1\otimes a_2\otimes a_3, \dot y_1\otimes a_2\otimes a_3 \rangle \\
&\enskip\quad + \lambda^2\langle  a_1\otimes \dot x_2\otimes a_3,  a_1\otimes \dot y_2\otimes a_3 \rangle + \lambda^2\langle a_1\otimes a_2\otimes \dot x_3, \dot a_1\otimes a_2\otimes \dot y_3 \rangle \nonumber \\
&=\dot\lambda\dot\mu + \lambda^2  \left( \langle \dot x_1,\dot y_1\rangle_{\mathbb{S}^{n_1-1}} +   \langle \dot x_2,\dot y_2\rangle_{\mathbb{S}^{n_2-1}}+\langle \dot x_3,\dot y_3\rangle_{\mathbb{S}^{n_3-1}}  \right)
,
\end{align*}	
where $\langle \cdot,\cdot \rangle_{\mathbb{S}^{n_i}}$ is the inner product on the sphere, which corresponds to the usual Euclidean inner product. We used orthogonality in the second equality and the fact that $a_i\in \mathbb{S}^{n_i-1}$ and hence has unit norm in the third equality. We can generalize these computation in the following lemma.
\begin{lem} \label{lemmametric} Let $p=(\lambda,a_1,\dots,a_d) \in Seg$ and take the tangent vectors\\ $(\dot\lambda,\dot x_1, \dots, \dot x_d)$ and $(\dot\mu,\dot y_1, \dots, \dot y_d) \in T_pSeg$. The metric induced by the ambient Euclidean space on $Seg$ is given by
\[ \langle (\dot\lambda,\dot x_1, \dots, \dot x_d),(\dot\mu,\dot y_1, \dots, \dot y_d)\rangle_p = \dot\lambda\dot\mu+\lambda^2\left( \sum_{i=1}^d\langle \dot x_i,\dot y_i \rangle_{\mathbb{S}^{n_i-1}} \right). \]
\end{lem}
Endowing the Segre manifold with this induced metric turns it into a Riemannian manifold. Note that this metric is closely related to the product metric of the usual metrics on $\mathbb{R}$ and $\mathbb{S}^{n_1-1}\times\mathbb{S}^{n_2-1}\times\mathbb{S}^{n_3-1}$. If $M$ and $N$ are two Riemannian manifolds, then one can endow $M\times N$ with a \textit{warped} product as follows: Given two tangent vectors $(v_1,w_1)$ and $(v_2,w_2)$ at a point $(p,q)\in M\times N$, one can define the inner product as 
\[ \langle (v_1,w_1),(v_2,w_2) \rangle _{(p,q)} = \langle v_1,v_2 \rangle_M + f^2(p)\langle w_1,w_2\rangle_N, \]
where $f: M\to\mathbb{R}$ is called the \textit{warping function}. The manifold endowed with this warped product is denoted by $M\times_f N$. In the case of $f$ being the constant function $1$, this reduces to a regular product manifold. Using this notation one can see that \[ Seg \cong \mathbb{R}_0\times_{Id_\mathbb{R}}(\mathbb{S}^{n_1-1}\times\mathbb{S}^{n_2-1}\times\mathbb{S}^{n_3-1}),\] 
hence it is a cone over a product of spheres. For more details on warped product we refer to \cite{oneill_semi-riemannian_1983}.

The remainder of this section serves as a proof of \cref{geo}. 
\begin{proof}[Proof of \cref{geo}]
Let us look at a path $ \gamma(t) = (\lambda(t),x_1(t),\dots,x_d(t)) 
\in \mathbb{R}\times \mathbb{S}^{n_1-1} \times \dots \times \mathbb{S}^{n_d-1}$ and assume it is the geodesic with given starting conditions $\left(
\lambda(0),x_1(0),\dots,x_d(0)\right)
$ and $(\dot\lambda(0), \dot{x_1}(0), \dots, \dot{x_d}(0))$.

Recall from the introduction of this section that the Segre manifold endowed with the induced metric is a warped product of $\mathbb{R}$ endowed with the standard metric and $\mathbb{S}^{n_1-1}\times \dots\times \mathbb{S}^{n_3-1}$ endowed with the product metric, using the identical warping function. This implies, as in \cite{oneill_semi-riemannian_1983}, that the path $\left( x_1(t), \dots , x_d(t) \right)$ is a \textit{pre-geodesic}, i.e., a geodesic that is not necessarily arc-length parametrized, in $\mathbb{S}^{n_1-1} \times \dots \times \mathbb{S}^{n_d-1}$. 

Now let us re-parametrize the path $\gamma$ in such a way that the angular path on the spheres has constant velocity $1$. Let us thus apply the change of variable $t \to t(s)$ with $t(0)=0$ such that we obtain the path  $\sigma(s) =  \left(\beta(t(s)),y_1(t(s)),\dots,y_d(t(s))
\right)$  with $\|(\dot y_1 , \dots, \dot y_d)\| =1$. It holds for all $1\leq i \leq d$ that
\[ ||\dot y_i(0)|| = ||\dot x_i(0)||t'(0) \quad \text{and} \quad \sum_{k=1}^d||\dot y_k||^2 =1. \]
From this we find that
\[ t'(0)  = \frac{1}{\sqrt{\sum_{k=1}^d||\dot y_k||^2}} =: \frac{1}{M}. \]
After this re-parametrization, the path $\left(y_1(t(s)), \dots, y_d(t(s))\right)$ is a proper geodesic, since it is an arc-length parametrization of a pre-geodesic. It is well known that this implies that the projection of this geodesic on each of the spheres results in a geodesic. Since geodesics on spheres are well-known, we can conclude that 
\[ y_i(s) = x_i(0)\cos\left(\frac{||\dot{x}_i(0)||}{M}s\right)+\frac{\dot x_i(0)}{||\dot x_i(0)||}\sin\left(\frac{||\dot x_i(0)||}{M}s\right), \]
for all $1\leq i \leq d$, where $M=\sqrt{\sum_{k=1}^d||\dot y_k||^2}$.

We now aim to find the path $\beta(t(s)): \mathbb{R} \to \mathbb{R}$. We will use the fact that geodesics are curves on a manifold that locally minimize distances. Recall that the length of $\sigma$ from $0$ to $T$ is 
\begin{equation}
\label{distance} \int_0^T \sqrt{\langle \dot\sigma(s),\dot\sigma(s) \rangle_{\sigma(s)}}\mathrm{d}s. 
\end{equation}
This shows that we aim to minimize \cref{distance} for all paths of the form $ \sigma(s)=\left(\beta(s), y_1(s),\dots, y_d(s)\right)$  with $||(\dot y_1,\dots, \dot y_d)|| =1$. Using the induced metric described in \cref{lemmametric}, \cref{distance} can be rewritten as
\begin{align*}
\int_0^T \sqrt{\langle \dot\sigma(s),\dot\sigma(s) \rangle_{\sigma(s)}}\ \mathrm{d}s &=\int_0^L \sqrt{\dot\beta\dot\beta + \beta^2  \left(\langle \dot y_1,\dot y_1\rangle_{\mathbb{S}^{n_1-1}}+\dots+\langle \dot y_d,\dot y_d\rangle_{\mathbb{S}^{n_d-1}}  \right)}\ \mathrm{d}s \\			
& =\int_0^L \sqrt{\dot\beta^2(s)+ \beta^2(s)}\ \mathrm{d}s.
\end{align*}
This minimization problem can be solved via the \textit{Euler--Lagrange equations} \cite{fox_introduction_1987}. When we define $L(s,\beta,\dot\beta):= \sqrt{\dot\beta^2(s)+ \beta^2(s)}$, we need to solve 
\[ \frac{\partial L}{\partial \beta} = \frac{\partial}{\partial s}\left(\frac{\partial L}{\partial \dot \beta} \right). \]
A straightforward computation gives 
\begin{align*}
&\frac{\partial L}{\partial \beta} = \frac{\partial}{\partial s}\left(\frac{\partial L}{\partial \dot \beta} \right) \\
&\Rightarrow\frac{\beta(s)}{\sqrt{\dot\beta^2(s)+ \beta^2(s)}} = \frac{\partial}{\partial s}\left(\frac{\dot\beta(s)}{\sqrt{\dot\beta^2(s)+ \beta^2(s)}}\right)\\
&\Rightarrow \frac{\beta(s)}{\sqrt{\dot\beta^2(s)+ \beta^2(s)}} = \frac{\ddot\beta(s)\sqrt{\dot\beta^2(s)+\beta^2(s)} - \frac{\dot\beta(s)\left(\dot\beta(s)\ddot\beta(s)+\beta(s)\dot\beta(s)\right)}{\sqrt{\dot\beta^2(s)+\beta^2(s)}}}{\dot\beta^2(s)+\beta^2(s)} \\
&\Rightarrow\beta(s)\left(\dot\beta^2(s)+\beta^2(s)\right) = \ddot\beta(s)\left(\dot\beta^2(s)+\beta^2(s)\right)- \dot\beta(s)\left(\dot\beta(s)\ddot\beta(s)+\beta(s)\dot\beta(s)\right),\end{align*}
which, after cancellations, reduces to
\[ \beta^2+2\dot\beta^2=\beta\ddot\beta.\]
One can check that the solution to this ODE is 
\[ \beta(s) = \frac{A}{\cos(s+B)}, \]
where $A$ and $B$ are constants to be determined by the starting conditions. Note that the starting conditions $\lambda(0)$ and $\dot\lambda(0)$ were given for the original path $\gamma$ and that $t(0)=0$. Hence we have
\[\begin{cases}
\begin{aligned}
\beta(0)&=\lambda(0)\\
\dot\beta(0)&=\dot\lambda(0)t'(0)
\end{aligned}
\end{cases}\Longrightarrow
\begin{cases}\begin{aligned}
\frac{A}{\cos(B)}&=\lambda(0) \\
\frac{A\sin(B)}{\cos^2(B)}&=\frac{\dot\lambda(0)}{M}
\end{aligned}\end{cases}
\Longrightarrow \begin{cases}\begin{aligned}
A&=\frac{\lambda(0)}{\sqrt{\left(\frac{\dot\lambda(0)}{\lambda(0)M}\right)^2+1}} \\
B&=\tan^{-1}\left(\frac{\dot\lambda(0)}{\lambda(0)M}\right).
\end{aligned}\end{cases}\]
We have thus found a pre-geodesic of the form
\begin{equation} \label{pregeo} 
\begin{split}
\beta(s) &= \frac{\lambda(0)}{\sqrt{\left(\frac{\dot\lambda(0)}{\lambda(0)M}\right)^2+1}}\frac{1}{\cos\left(s+\tan^{-1}\left(\frac{\dot\lambda(0)}{\lambda(0)M}\right)\right)}, \\
y_i(s)&=y_i(0)\cos\left(\frac{||\dot{y_i}(0)||}{M}s\right)+\frac{\dot y_i(0)}{||\dot y_i(0)||}\sin\left(\frac{||\dot{y_i}(0)||}{M}s\right).
\end{split}
\end{equation}
All that remains is computing an arc-length parametrization of this path. Denote this parametrization by $s  \to s(u)$, with the assumption that $s(0)=0$, such that the path $\chi(u):= \sigma(s(u))$ has length $1$. Using the abbreviation $P:=\frac{\dot\lambda(0)}{\lambda(0)M}$, we can compute
\begin{align*}
1&= ||\dot\chi(u)|| \\
&= ||\dot\sigma(s(u)))||s'(u) \\
&= s'(u)\sqrt{\frac{\lambda^2(0)}{(P^2+1)}\frac{\sin^2(s+\tan^{-1}(P))}{\cos^4(s+\tan^{-1}(P))}+\frac{\lambda^2(0)}{(P^2+1)}\frac{1}{\cos^2(s+\tan^{-1}(P))}} .
\end{align*}
The solution of this ODE is straightforward via separation of variables:
\begin{align*}
\mathrm{d}u&=\frac{\lambda(0)}{\sqrt{P^2+1}}\frac{1}{\cos^2(s+\tan^{-1}(P))}\mathrm{d}s\quad; \\
\text{so that}\quad u+C &= \frac{\lambda(0)}{\sqrt{P^2+1}}\tan(s+\tan^{-1}(P)).
\end{align*}
Using that $s(0)=0$, we can compute $C$ to be $\frac{\lambda(0)P}{\sqrt{1+P^2}}$. We can then rewrite to conclude that
\[ s= \tan^{-1}\left(\frac{\sqrt{P^2+1}}{\lambda(0)}u+P\right)-\tan^{-1}(P). \]
We can apply this change of variable to obtain the arc-length parametrized geodesic, which takes the following form after simplifying and renaming $u$ to $t$:
\begin{align*}
\lambda(t)&=\sqrt{t^2+\frac{2\lambda(0)P}{\sqrt{P^2+1}}t+\frac{\lambda^2(0)}{P^2+1}},\\
x_i(t)&=x_i(0)\cos\left(\frac{||\dot{x_i}(0)||}{M}\left(\tan^{-1}\left(\frac{\sqrt{P^2+1}}{\lambda(0)}t+P\right)-\tan^{-1}\left(P\right)\right)\right) \\
&+\frac{\dot x_1(0)}{||\dot x_1(0)||}\sin\left(\frac{||\dot{x_i}(0)||}{M}\left(\tan^{-1}\left(\frac{\sqrt{P^2+1}}{\lambda(0)}t+P\right)-\tan^{-1}\left(P\right)\right)\right),
\end{align*}
for all $1\leq i\leq d$, with $M=\sqrt{\sum_{k=1}^d||\dot{x}_k||^2}$ and $P=\frac{\dot\lambda(0)}{\lambda(0)M}$. Using the abbreviations as in \cref{geo}, we obtain the result. 
\end{proof}
\section{Riemannian conjugate gradient algorithm for tensor decomposition}
\label{sectionalgo}
In the previous section we had an in-depth look at the geodesics of the Segre manifold. The aim of this section is to bring together all items needed to write a Riemannian conjugate gradient algorithm to find the best rank $r$ decomposition of a given tensor.

A conjugate gradient algorithm is an improvement of the gradient descent algorithm, since it remembers the previous search direction. A gradient descent method on an embedded manifold has a clear geometric interpretation. One starts with a gradient flow in the ambient space pointing towards a point. One then projects this flow onto the submanifold and aims to follow this projected flow.

Since we are optimizing over a manifold, we need to compute certain geometric properties of the manifold. We will use the geodesics to move over the manifold in a given direction. Additionally we will use parallel transport to move the previous search direction to the current iteration point. The pseudo-code below is a basic layout of a conjugate gradient algorithm over a manifold; for more details see	 \cite{absil_optimization_2009}.

\begin{algorithm}[H]
\SetAlgoLined	

\KwData{Manifold $M$, function $f:M\to \mathbb{R}$}
\KwResult{Minimum of function $M$}
Pick a starting point $p\in M$\;
\While{Stopping conditions not met}{
Compute Riemannian gradient $v \in T_pM$\;
Move previous gradient to $T_pM$ with parallel transport to obtain $w$\;
Compute $v_{new} = -(v + \beta w)$\;
Use line-search algorithm to find step-size $\alpha$\;
Use retraction to move from $p$ in direction of $v_{new}$ for time $\alpha$ to obtain $p_{new}$\;
}
\caption{Riemannian conjugate gradient algorithm}
\end{algorithm}

There are different methods for choosing the parameter $\beta$ in the above algorithm. Our algorithm uses a Hestenes-Stiefel rule.  For details we again refer to \cite{absil_optimization_2009}.

\subsection{Objective function}

Recall that the goal is to find the best rank-$r$ approximation of a given tensor $T$. We aim to solve the following least squares optimization problem:
\[\min_{(T_1,\dots,T_r)\in Seg^{\times r}} \frac{1}{2}\left\|T-\sum_{i=1}^rT_i\right\|^2_{\text{F}}. \]
As explained before, we optimize over the product manifold of $r$ copies of $Seg$, denoted by $Seg^{\times r}$, as in \cite{breiding_riemannian_2018}. In the remainder of this section, we often use $M$.  

\subsection{Riemannian gradient}

To compute the search direction from a point $p\in M$, we need to compute the Riemannian gradient of the objective function at $p$. The Riemannian gradient of $f:M \to \mathbb{R}$ is defined as the unique vector field $\nabla f$ such that for all $p\in M$ and $v \in T_pM$:
\[ v(f) = \langle v,\nabla f\rangle.  \]
Recall that $Seg$ is an embedded manifold in the ambient Euclidean space $\mathbb{R}^{n_1\times\dots\times n_d}$. It is shown in \cite{boumal_introduction_2020}, that the Riemannian gradient at a point $p$ on the Segre is the projection of the gradient of the extended map $\tilde f$ on the ambient space to $T_pM$. Thus we can conclude that 
\[ \nabla f_{(T_1,\dots,T_r)} = P_{(T_1,\dots,T_r)} (\nabla\tilde f_{(T_1,\dots,T_r)}) = P_{(T_1,\dots,T_r)} \left(T-\sum_{i=1}^r T_i\right), \]
where $P_p$ is the projection onto $T_pM$. This gradient can be computed efficiently as in \cite{vannieuwenhoven_computing_2015}.

\subsection{Geodesics}

In a Riemannian conjugate gradient algorithm we often need to move from a point $p\in M$ in a certain direction $v\in T_pM$. We will use the geometrically optimal curves
\begin{align*}
R_{(p_1,\dots,p_r)}: T(Seg^{\times r}) &\to Seg^{\times r} \\
(v_1,\dots,v_r) &\mapsto (\gamma_1(1),\dots,\gamma_r(1))
\end{align*}
where $\gamma_i:T_{p_i}Seg \to Seg$ is the geodesic starting at $p_i$ in the direction $v_i$. This map is known as the $\emph{exponential map}$. 

\subsection{Parallel transport}

Additionally, a Riemannian conjugate gradient algorithm requires a method to move tangent vectors from one iteration point to the next. The geometric tool corresponding to this is a connection $\nabla$ on the manifold $M$. Note that the Segre manifold is an embedded manifold in the ambient Euclidean space of all tensors of the fixed size. It is described, for example in \cite{boumal_introduction_2020}, that the parallel transport is simply the projection of the ambient parallel transport onto the tangent spaces. If we call the ambient connection $\overline{\nabla}$ and take two vector fields $U,V$ on $M$, we find
\[ \nabla_U V = P_p\left(\overline{\nabla}_U \overline{V}\right) \]
where $\overline{V}$ is a smooth extension of the vector field to the ambient Euclidean space and $P_p$ is the orthogonal projection onto $T_pSeg$. 

\subsection{Line search}
\label{linesearch}
Once we have determined the direction in which we move from a certain iteration point, we need a line search algorithm to dictate the step size $\alpha$. A commonly used line search algorithm to determine a good step size $\alpha$ is a so-called backtracking algorithm. One typically, as is mentioned in \cite{boumal_introduction_2020}, starts with an initial guess for $\alpha$ and iteratively decreases it with a factor (often $0.5$) until the \textit{Armijo--Goldstein} condition is satisfied \cite{boumal_introduction_2020}. This condition checks if the decrease of the objective function using a certain step size is at least what we would expect given the step size and the norm of the gradient.

\begin{figure}[h]
\includegraphics[scale=0.3]{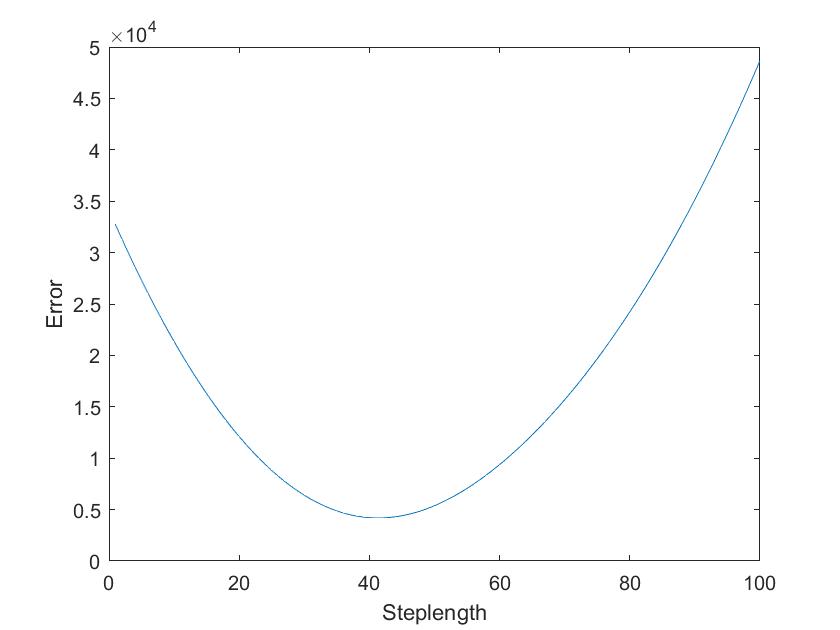}
\caption{Error for different step-sizes}
\label{linesearchfig}
\end{figure}

One disadvantage of this method is the unknown amount of iterations and therefore the uncontrolled amount of function evaluations. After investigating the particular problem of tensor approximation, it was noted that the evaluation of the objective function behaved nicely with respect to the step size, as can be seen in \cref{linesearchfig}.
Based on these observations we propose a simple formula using quadratic interpolation. One already knows the current function value corresponding to a step size $\alpha=0$, call this $f_0$. Computing the function values corresponding to $\alpha=1$ and $\alpha=2$, called respectively $f_1$ and $f_2$ one can easily calculate the interpolating parabola as
\[ f(\alpha) = \frac{(\alpha-1)(\alpha-2)}{2}f_0-\alpha(\alpha-2)f_1+\frac{\alpha(\alpha-1)}{2}f_2.\]
This parabola obtains its minimum at step size
\[ \alpha = \frac{3f_0-4f_1+f_2}{2f_0-4f_1+2f_2}. \]

Since the error at step size $0$ has to be calculated to determine the error at the current iteration point, this method requires only $2$ additional function evaluations. This clear bound on the number of function evaluations is the main advantage of this method of determining $\alpha$. In the code for this line search algorithm we still check if the $\alpha$ that is determined by the quadratic interpolation satisfies the Armijo--Goldstein condition. If not, we return to the default backtracking line search algorithm. 

\section{Experimental results}
\label{sectionresults}
In this section we apply the Riemannian conjugate gradient algorithm to two different topics, namely recommender systems for movies and fluorescence spectroscopy. 

The Riemannian conjugate gradient algorithm using geodesics on the Segre manifold is implemented in Matlab. The toolbox Manopt \cite{boumal_manopt_2013} was used for the general framework of optimization over manifolds. In particular, we used the \texttt{conjugategradient} method with the default options. Note that instead of parametrizing $Seg^{\times r}$ as $\left(\mathbb{R}\times \mathbb{S}^{n_1-1}\times \mathbb{S}^{n_2-1} \times \mathbb{S}^{n_3-1}\right)^{\times r}$, we worked with  \[ Seg^{\times r} \cong \mathbb{R}^r \times \mathcal{OB}(r,n_1) \times \mathcal{OB}(r,n_2)\times \mathcal{OB}(r,n_3), \]
where $\mathcal{OB}(r,n_1)$ is the \textit{oblique manifold} consisting of all $n_1\times r$ matrices with unit norm column vectors. Its geometry is equivalent to $\left(\mathbb{S}^{n_1-1}\right)^{\times r}$, but its implementation is more efficient.  Additionally, tools from Tensorlab \cite{vervliet_tensorlab_2016} were used for efficient tensor calculations. 

The experiments were run on an HP Elitebook. It has an AMD Ryzen $7$ PRO $2.30$ GHz processor and $32$ GB of RAM. 

\subsection{Recommender system for MovieLens data}

In this well known application we will apply the algorithm to answer the problem of tensor completion. In the field of tensor completion, one aims to complete a given incomplete tensor, assuming an underlying low rank structure \cite{song_tensor_2019}. The dataset that will be used is the 1M MovieLens dataset \cite{harper_movielens_2015}. This database contains $1\, 000\, 209$ ratings (ranging from $1$ to $5$ stars) from $6040$ users of $3900$ movies. In order to increase the density, the data have already been filtered by MovieLens in such a way that every user has rated at least $20$ movies. The density of the data, i.e., the percentage of the known elements, is approximately $4.2\%$.

The aim of a recommender system is to predict the match between an arbitrary user and movie as accurately as possible. We are thus interested in a recommendation map:
\[ f_r : \text{Users} \times \text{Movies} \to \text{RelevanceScore}. \]
As is done in \cite{frolov_fifty_2016}, we interpret the dataset as a $6040 \times 3900 \times 5$ tensor $T$. The first two modes of $T$ correspond to the users and the movies, while the last mode corresponds to the ratings $\{ 1,2,3,4,5\}$. So if user $i$ gives movie $j$ a rating of $3$ stars, we will have $T(i,j,3)=1$ and $T(i,j,k)=0$ if $k \neq 3$. The aim of this process is instead of generating an expected rating for a given (user,movie)-pair, we can approximate a rating distribution. In other words, we aim to create a recommender map that looks like
\[ f_r: \text{Users}\times \text{Movies} \times \text{Ratings} \to \text{RelevanceScore}. \]
For more details and an overview of tensor methods being used in recommender systems, we refer to \cite{frolov_tensor_2017}. As is usual, we assume the data can be well-approximated by a low-rank tensor decomposition. We aim to find the best rank-$r$ approximation for the tensor $T$, of course only looking at the known indices. The tensor can then be completed by this model, thusly recommending ratings for unknown user-movie pairs. Before showing the experimental results, we outline some specifics of the used algorithm for this setting.

We start by discussing the objective function. As mentioned above we aim to minimize the distance between a rank $r$ tensor and $T$. To compute this distance, we only look at the set of given indices, which we call $I$. Note that intuitively we expect entries of the tensor to be values in the interval $[0,1]$. To promote this behaviour, we add a penalty function that punishes each entry of the tensor outside of this interval. Concretely, we define the objective function as
\[ f(A) = \frac{1}{2}\left\|T-A\right\|^2_{\text{F,I}}+\lambda\sum_{i=1,j=1,k=1}^{n_1,n_2,n_3}\left(A(i,j,k)^2-A(i,j,k)\right)^9, \]
where $A$ is the low rank approximation and $\left\|\cdot\right\|_{\text{F,I}}$ is the Frobenius norm obtained by only summing over the indices in $I$. The parameter $\lambda$, which is taken equal to $1$ in the remainder of this section, is used to give weight to the penalty term. The ninth power in the penalty term is used to push the negative part of the graph of $x^2-x$ to zero in a smooth way. The gradient of the objective function was adjusted accordingly. The addition of the penalty function removes the nice effect observed in \cref{linesearch}, hence we used the standard backtracking algorithm provided by Manopt. We used the \texttt{rand} fucntion by Manopt to randomly choose a starting point on the product manifold.

Secondly we need a way to recommend a rating $R$ for a user-movie pair $(i,j)$, using our five predicted values $A(i,j,k)$. We considered four possible methods:

\begin{enumerate}
\item Defining $R$ as the weighted average, so $R=\sum_{k=1}^5 k*A(i,j,k)$.
\item First rescaling the scores to have sum $1$, then proceeding as in the first option.
\item Setting negative values equal to zero and setting values bigger than one equal to one, then proceeding as in the first option.
\item Recommending a rating of $k$ stars if the highest value is achieved is $A(i,j,k)$.  
\end{enumerate}

\begin{figure}[h!]
\includegraphics[scale=0.3]{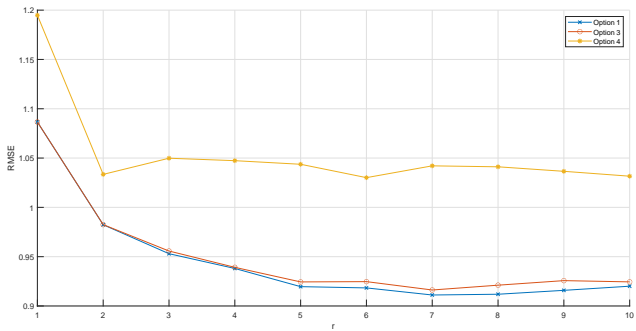}
\caption{RMSE for different methods of predicting ratings.}\label{rmse}
\end{figure}

As is customary, the full dataset was randomly split $80\%/20\%$ into a training set and a test set. Only the training set was used to train the algorithm. Firstly, one aims to find the optimal $r$ for approaching the tensor $T$ by a rank $r$ tensor. Choosing a small value for $r$ risks not being able to properly explain the dataset, while choosing a large value for $r$ risks overfitting the data for the training set and obtaining worse recommendations for the test set. For each rank $r$ we make a new random $80\%/20\%$ split of the training set. We run the algorithm to approximate the $80\%$ data and compare the prediction to the $20\%$ data we withheld. We used the notion of the \textit{root mean square error} (RMSE) to compare the predictions to the real values. RMSE is a well-established measure of difference between observed values and values predicted by a model and is given by
\[ RMSE(x,y) = \sqrt{\frac{\sum_{i=1}^n (y_i-x_i)^2}{n}} \]
We stopped the algorithm when the norm of the gradient dropped below $10^{-1}$. In all of our runs, this corresponds to a relative norm of about $10^{-5}$. One can find the results in \cref{rmse}. Note that option $2$ is not shown, since it was not at all competitive with RMSE up to $200\%$ larger than the other options.  

It is clear that option $1$ performs best in this setting, with option $2$ being the only other competitive option. Already for reasonably low ranks, such as $r=7$, we achieve values as low as $0.9017$. These values are in the range $[0.875,0.964]$ of  several different algorithms; see \cite{papadakis_scor_2017}. 

Running the algorithm for $r=7$ takes around $250$ iterations in $11$ minutes. Note however that this is not the relevant time. One is for example interested in quickly giving the top recommendations for a given user. Taking $100$ arbitrary users and providing each of them with the $10$ films best suited for them took $0.010$ seconds.

\subsection{Fluorescence spectroscopy}

As a final application we look at a dataset generated by fluorescence spectroscopy. The experiments aim to identify the concentration of several fluorophores in different mixtures. This is done by exciting them with light of different wavelengths and intercepting the light emitted by the mixtures.

Currently, experiments in fluorescence spectroscopy are executed by exciting the mixtures with a range of wavelengths and intercepting a range of wavelengths. It is of interest to determine how much data is really required to identify the fluophores in the different mixtures. Concretely we ask ourselves how incomplete the tensor with data can be in order to still be able to recover the full dataset. This can be useful to develop new hardware for fluorescence spectroscopy. By only exciting the mixtures with light of randomly chosen wavelengths and using a filter to only intercept light at certain random wavelengths, we anticipate that one could drastically reduce the experimental duration. We investigate this on two different datasets: one smaller dataset with a small relative error and one bigger dataset with a larger relative error.

It is known that one can expect a dataset obtained by fluorescence spectroscopy on mixtures of $r$ different fluophores to be explained by a model of rank $r$ \cite{appellof_three-dimensional_1983}. We will check this experimentally in both datasets. We aim to experimentally determine how incomplete the tensor of data can become such that it can still be  recovered. We will remove all but $\phi\cdot\dim(\sigma_r(Seg))$ elements of $T$, where $\phi$ is a variable and $\sigma_r(Seg)$ is the $r$-th secant variety of $Seg$. Note that for order three tensors of size $n_1\times n_2\times n_3$ one has $\dim(\sigma_r(Seg))\leq r(n_1+n_2+n_3-2)$, with equality generally expected to hold for sufficiently small ranks $r$ \cite{abo_induction_2009}. Results will be quantified using both the relative approximation error and the core consistency of the approximation, which is a constant derived from a procedure called \textit{Corcondia}\cite{bro_new_2003}. Corcondia is often used for deriving the amount of latent factors in a model, by computing a percentage that decreases when a model is being overfitted. For details on this method, together with an efficient implementation we refer to \cite{papalexakis_fast_2015}.   

\begin{figure}[h!]
\includegraphics[scale=0.4]{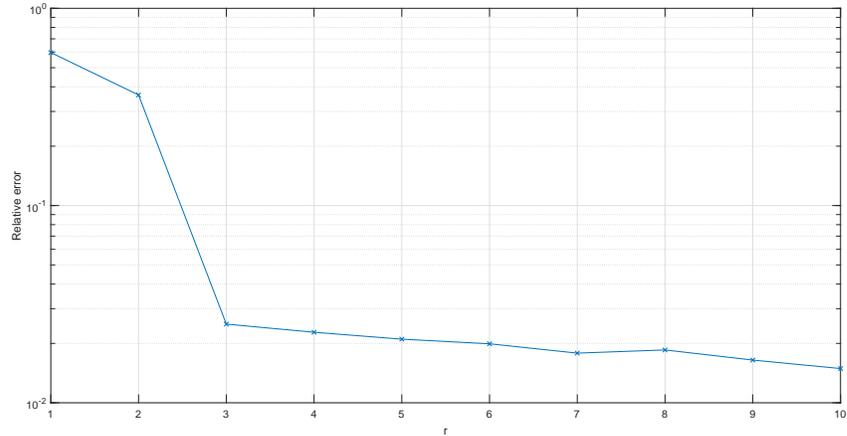}
\caption{Relative errors per rank}
\label{relAmino}
\end{figure}

\begin{figure}[h]
\includegraphics[scale=0.3]{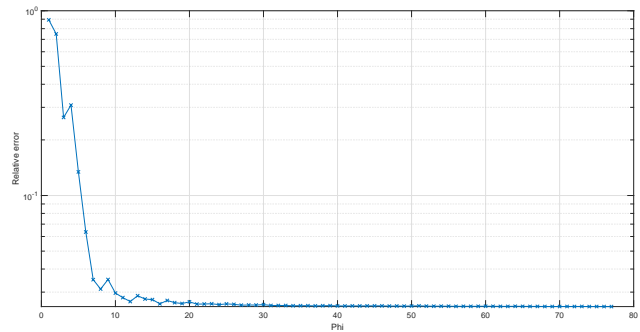}
\caption{Relative error of the approximation obtained using only $\phi\cdot\dim(\sigma_r(Segre))$ elements.}
\label{aminophi}
\end{figure}

\begin{figure}[h]
\includegraphics[scale=0.3]{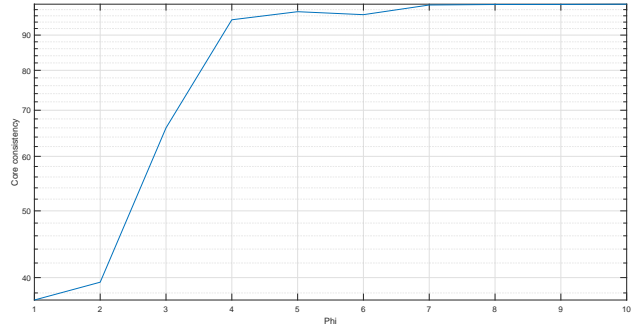}
\caption{Core consistency of the approximation obtained using only $\phi\cdot\dim(\sigma_r(Segre))$ elements.}
\label{CCAmino}
\end{figure}

\begin{figure}[h]
\begin{subfigure}[b]{0.475\textwidth}
\centering
\includegraphics[width=\textwidth]{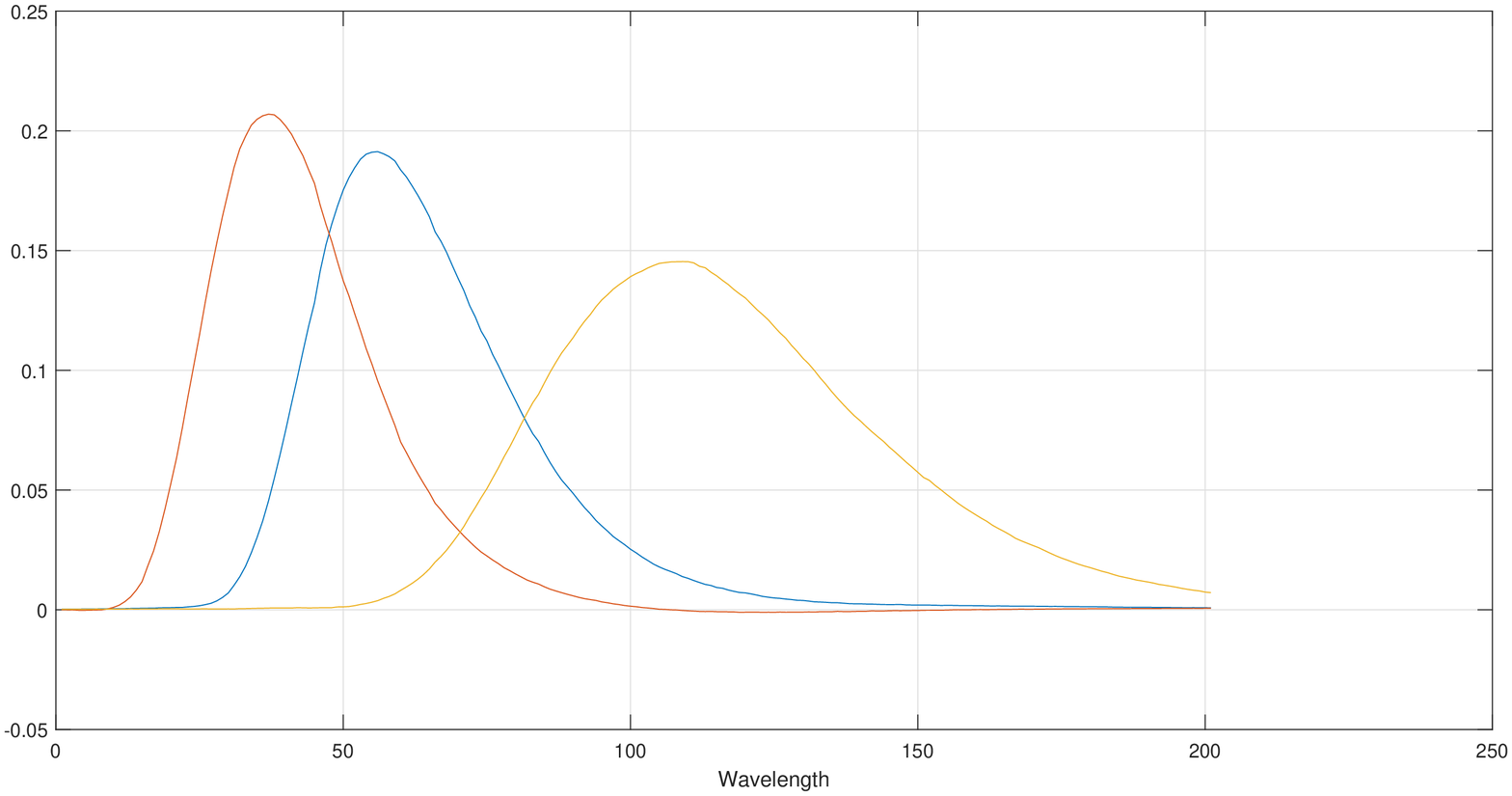}
\caption{Full tensor, relative error of $0.0251$}
\end{subfigure}
\hfill
\begin{subfigure}[b]{0.475\textwidth}
\centering
\includegraphics[width=\textwidth]{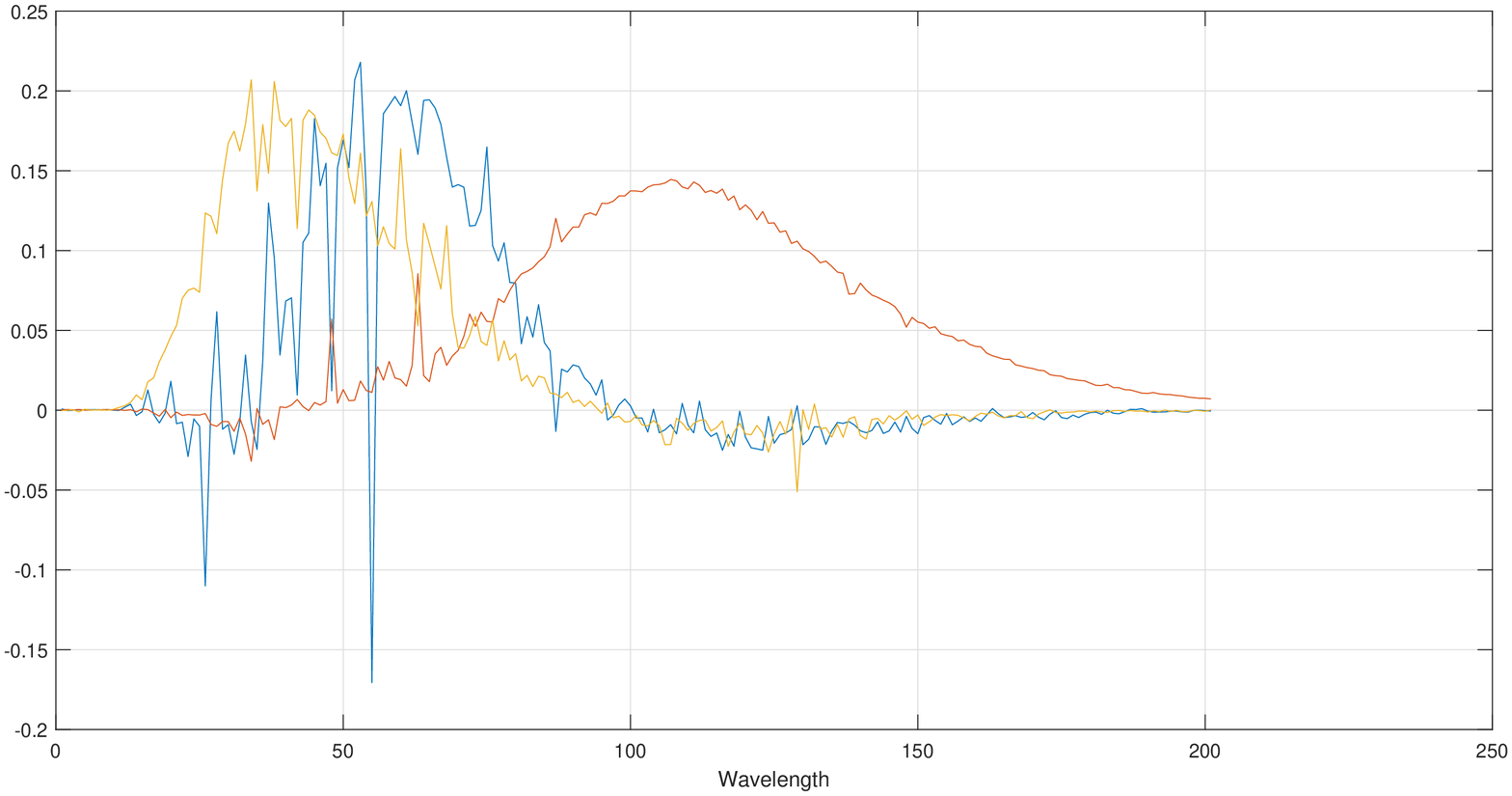}
\caption{$\phi=5$, relative error of $0.1342$}
\end{subfigure}
\vskip\baselineskip
\begin{subfigure}[b]{0.475\textwidth}
\centering
\includegraphics[width=\textwidth]{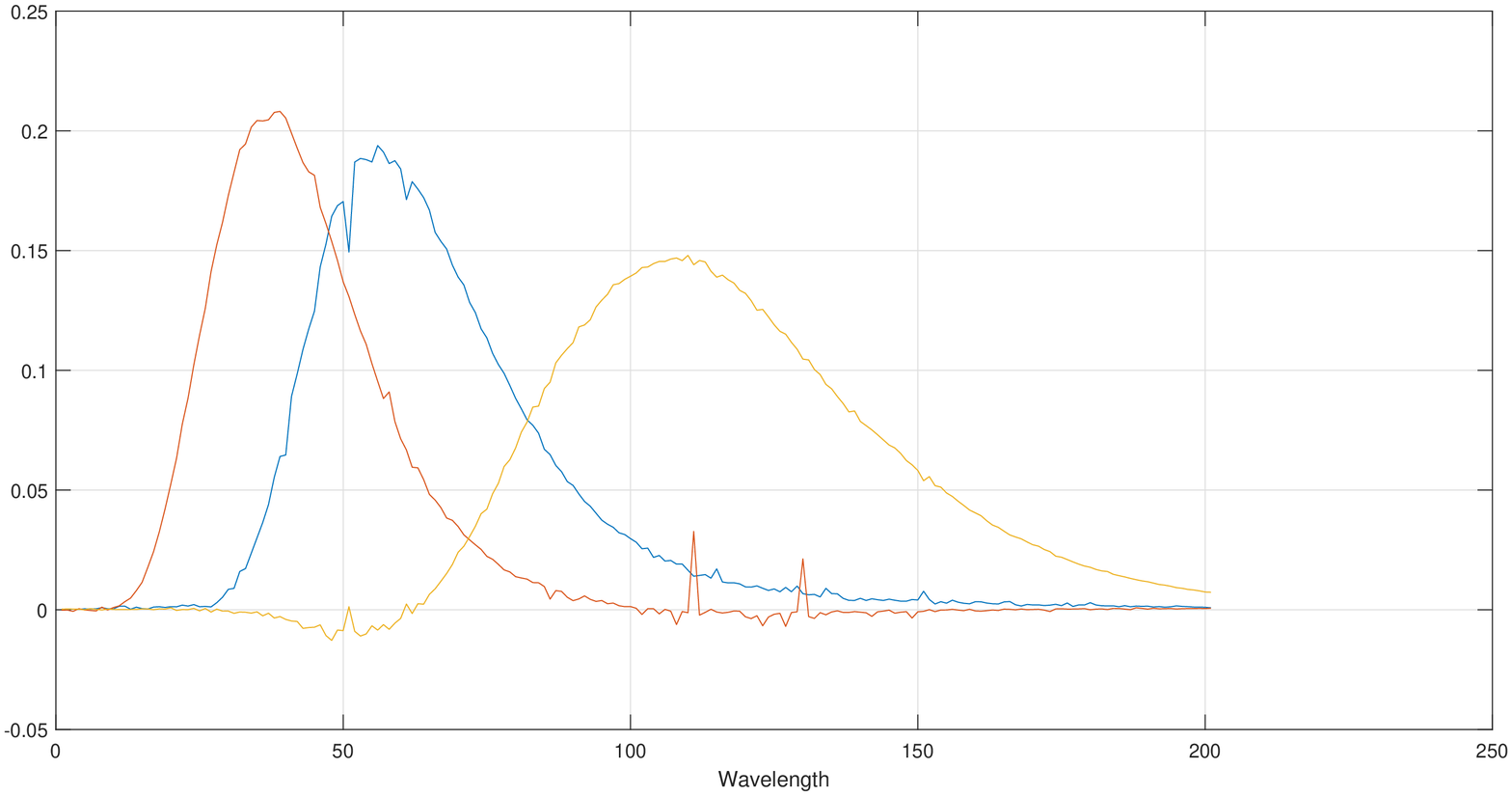}
\caption{$\phi=7$, relative error of $0.0350$}
\end{subfigure}
\hfill
\begin{subfigure}[b]{0.475\textwidth}
\centering
\includegraphics[width=\textwidth]{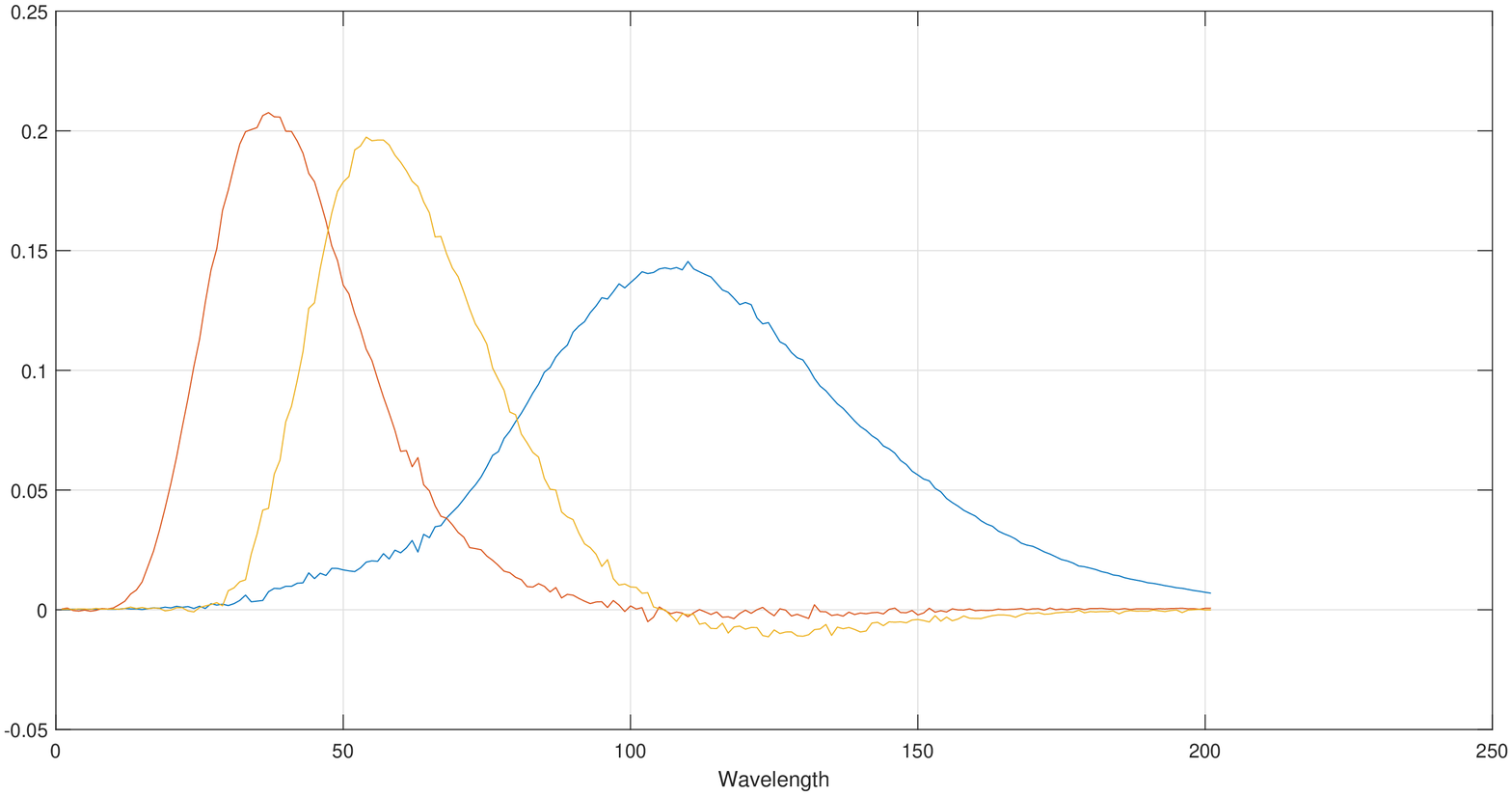}
\caption{$\phi=10$, relative error of $0.0296$}
\end{subfigure}
\caption{Emission patterns of the rank 3 decomposition}
\label{emission}
\end{figure}

The dataset \cite{bro_parafac_1997} consists of five mixtures, each containing different concentrations of three amino acid fluorescents: tyrosine, tryptophan and phenylalanine. These mixtures are excited by light of wavelengths between $240$ nm and $300$ nm in intervals of $1$ nm. The emitted light is intercepted at wavelengths between $250$ nm and $450$ nm in intervals of $1$ nm. The data can therefore be collected in a tensor $T$ of size $5\times 201\times 61$. \Cref{relAmino} shows the relative errors obtained by approximating $T$ with a rank $r$ tensor. Note that we stopped the algorithm when the norm of the gradient was less then $10$. In all of our runs, this corresponds to a relative norm of about $10^{-4}$. 
It is clear that, as expected, a model of rank $3$ is sufficient to explain the dataset, which corresponds to a relative error of around $0.025$.   Note that in this setting $\dim(\sigma_r(Seg)) = r(n_1+n_2+n_3-2)= 795$, which is $1.3 \%$ of the size of the entire dataset. \Cref{aminophi} shows the relative error obtained by only using $\phi\cdot\dim(\sigma_r(Seg))$ randomly chosen elements of $T$ in function of $\phi$. As a further study for values of $\phi \in [1,10]$, \cref{CCAmino} shows the core consistency of the approximation using only $\phi\cdot\dim(\sigma_r(Seg))$ elements.
\Cref{aminophi,CCAmino} suggest that for $\phi=7$ we already have an accurate approximation of $T$. One can check this by retrieving the emission patterns of the amino acids from the decompositions. \Cref{emission} shows that the emission patterns are well approximated from the incomplete tensor. This shows that in this setting one only needs around $9\%$ of all the data collected by fluorescence spectroscopy. Because of the relatively small dataset, this corresponds to a decrease in the runtime of the computer algorithm from $1.3$ seconds to $0.5$ seconds. However, since only about $9\%$ of the data has to be collected, there is a possibility to significantly lower the acquisition time of the spectroscopic images by developing new hardware.  

\section*{Conclusions}

We proposed a Riemannian conjugate gradient algorithm for tensor rank approximation using the geometry of the manifold over which one optimizes. The Riemannian geometry of the manifold $Seg$ of rank $1$ tensors was studied. The main original contribution is the computation of the geodesics of $Seg$ and applying them as retractions in the Riemannian conjugate gradient algorithm. We applied the algorithm in two different settings with incomplete tensors. First, we showed how to use the algorithm to build a recommender system for movie ratings. Second, we determined experimentally that one only needs a fraction of the data in order to recover a common fluorescence spectroscopy dataset up to measurement errors. This observation could be exploited to create new hardware for fluorescence spectroscopy that reduces the experimental duration.

\bibliographystyle{amsalpha}
\bibliography{referencesCG}

\end{document}